\documentclass[twosided,12pt,reqno,]{article}
\usepackage{amssymb}

\usepackage{amsthm}
\usepackage{amsmath}
\usepackage{a4wide}
\usepackage[all]{xy}
\usepackage{fancyhdr}

\newtheorem{theorem}{Theorem}[section]
\newtheorem{proposition}[theorem]{Proposition}
\newtheorem{corollary}[theorem]{Corollary}
\newtheorem{lemma}[theorem]{Lemma}
\theoremstyle{definition}
\newtheorem*{notation}{Notation}
\newtheorem*{Beweis}{Proof}
\newtheorem{definition}[theorem]{Definition}
\newtheorem{punto}[theorem]{}
\theoremstyle{remark}
\newtheorem{remark}[theorem]{Remark}
\newtheorem{ex}[theorem]{Example}
\newtheorem{exs}[theorem]{Examples}
\newtheorem{c-ex}[theorem]{Counterexample}
\newtheorem{c-exs}[theorem]{Counterexamples}
\newtheorem{remarks}[theorem]{Remarks}

\swapnumbers
\CompileMatrices
\input{tcilatex}

\begin{document}
\title{A Zariski Topology for Modules}
\author{{\bf Jawad Abuhlail}\footnote{The author would like to acknowledge the support provided by the Deanship of Scientific Research (DSR) at King Fahd University of Petroleum $\&$ Minerals (KFUPM) for funding this work through project No. FT08009.}\\ Department of Mathematics and Statistics\\ King Fahd University of Petroleum $\&$ Minerals \\ abuhlail@kfupm.edu.sa}
\date{\today}
\maketitle

\begin{abstract}
Given a duo module $M$ over an associative (not necessarily
commutative) ring $R,$ a Zariski topology is defined on the spectrum
$\mathrm{Spec}^{\mathrm{fp}}(M)$ of {\it fully prime} $R$-submodules of $M$.
We investigate, in particular, the interplay between the properties of
this space and the algebraic properties of the module under consideration.
\end{abstract}

\section{Introduction}

\qquad The interplay between the properties of a given ring and the topological
properties of the Zariski topology defined on its prime spectrum has been studied intensively
in the literature (e.g. \cite{LY2006, ST2010, ZTW2006}). On the other hand, many papers considered
the so called \emph{top modules,} i.e. modules whose spectrum of \emph{prime submodules} attains a
Zariski topology (e.g. \cite{Lu1984, Lu1999, MMS1997, MMS1998, Zha2006-a, Zha2006-b}).

On the other hand, different notions of primeness for modules were introduced and investigated
in the literature (e.g. \cite{Dau1978, Wis1996, RRRF-AS2002, RRW2005, Wij2006, Abu2006, WW2009}).
In this paper, we consider a notions that was not dealt with, from the topological point of view, so far.
Given a \emph{duo} module $M$ over an associative ring $R,$ we introduce and topologize
the spectrum of {\it fully prime} $R$-submodules of $M$. Motivated by results on the Zariski topology on the spectrum of prime ideals of a
commutative ring (e.g. \cite{Bou1998, AM1969}), we investigate the
interplay between the topological properties of the obtained space and the module under consideration. The result in this article are, in some sense, dual to those in \cite{Abu2008} in which the author investigated a Zariski-like topology for comodules over corings.

After this short introductory section, we introduce in Section 2 some
preliminaries. In particular, we recall some properties and notions from
Module Theory that will be needed in the sequel. In Section 3, we introduce
and investigate a Zariski topology on the spectrum $\mathrm{Spec}^{%
\mathrm{fp}}(M)$ of proper submodules that are \emph{fully prime} in $M$ (e.g. \cite{RRRF-AS2002}). In particular, and assuming suitable conditions for each result, we investigate when this space is Noetherian (Theorem \ref{noeth}), irreducible (Corollary \ref{rad-irred}), ultraconnected (Proposition \ref{ultra}), compact (Theorem \ref{fp-Lindelof}), connected (Theorem \ref{local}), $T_1$ (Proposition \ref{Frecht}) or $T_2$ (Theorem \ref{f.p.-discrete}). We end with applications to rings. In particular, we provide new characterizations of commutative $0$-dimensional semilocal rings and commutative semisimple rings (Corollary \ref{0-dim}) in terms of the so-called complete max-property which we introduce for modules (rings) in \ref{max-prop}.

\section{Preliminaries}

\qquad In this section, we fix some notation and recall some definitions and
basic results. For topological notions, the reader might consult any book in
General Topology (e.g. \cite{Bou1966}). For any undefined terminology, the
reader is referred to \cite{Wis1991} and \cite{AF1974}.

Throughout, $R$ is an associative (not necessarily commutative) ring
with $1_{R}\neq 0_{R}$ and $M$  is a non-zero unital left $R$-module. An ideal will mean a {\it two-sided ideal}, unless otherwise explicitly mentioned. By an $R$-module we mean a {\it left} $R$-module, unless explicitly otherwise mentioned. We set $S:=\mathrm{End}(_{R}M)^{op}$
(the ring of $R$-linear endomorphisms of $M$ with multiplication given by the opposite
composition of maps $(m)(gf)=(g\circ^{op}f)(m)=f(g(m))$ and consider $M$ as an $(R,S)$-bimodule in the
canonical way. We write $L\leq _{R}M$ ($L\lvertneqq _{R}M$) to indicate
that $L$ is a (proper)\ $R$-submodule of $M$. We say that $_R M$ is {\it distributive}, iff for any $L, L_1, L_2 \leq_R M$ we have $L \cap (L_1 + L_2) = (L \cap L_1) + (L \cap L_2)$. For non-empty subsets $L \subseteq M$ and $I \subseteq R$ we set $$(L :_R M):=\{r \in R | \, r M \subseteq L\} \, \text{   and   } \,  (L :_M I):=\{m \in M | \, I m \subseteq L\}.$$ Moreover, we set $$ {\rm Gen}(M):= \{ _R N \mid N \text{ is } M\text{-generated}\}; $$ $${\rm Cogen}(M):= \{ _R N \mid N \text{ is } M\text{-cogenerated}\}.$$

\begin{definition}
We say $L\leq _{R}M$ is \emph{fully invariant} or \emph{characteristic},
iff $f(L)\subseteq L$ for every $f\in S$ (equivalently,
iff $L\leq M$ is an $(R,S)$-subbimodule). In this case, we write $L\leq
_{R}^{\mathrm{f.i.}}M.$ We call $_{R}M$ \emph{duo} or \emph{invariant}, iff every $R$-submodule of $M$ is fully invariant.
\end{definition}

Recall that the ring $R$ is said to be {\it left duo} ({\it right duo}), iff every left (right) ideal is two-sided and to be {\it left quasi-duo}  ({\it right quasi-duo}) iff every maximal left (right) ideal of $R$ is two-sided. Moreover, $R$ is said to be ({\it quasi-}) {\it duo}, iff $R$ is left and right (quasi-) duo.

\begin{punto}\label{ex-duo}

Examples of duo modules are:

\begin{enumerate}

\item uniserial Artinian modules over commutative rings \cite{OHS2006}.

\item {\it multiplication modules}: these are $R$-modules $M$ such that every $L \leq_R M$ is of the form $L=IM$ for some ideal
$I$ of $R$, equivalently $L = (L:_R M)M$. It is obvious that every multiplication module is duo. Multiplication modules over commutative rings have been studied intensively in the literature (e.g. \cite{AS2004, PC1995, Smi1994}). By results of Tuganbaev \cite{Tug2003} (see also \cite{Tug2004}), we have the following generalizations to rings {\it close to be commutative}:

\begin{enumerate}

\item If $R$ is a left duo ring,  then the following classes of left $R$-modules are multiplication: cyclic modules; finitely generated distributive modules; left ideals which are generated by idempotents.

\item If $R$ is left quasi-duo, then every finitely generated semisimple distributive left $R$-module is multiplication.

\item If $R$ is duo,  then invertible ideals are finitely generated multiplication left (and right) $R$-modules.

\item If $R$ is a ring with commutative multiplication of left ideals (e.g. $R$ is {\it strongly regular}), then all projective left ideals of $R$ are multiplication and every finitely generated $R$-module $M$ for which $M_P$ is a cyclic $R_P$-module for every maximal left ideal of $R$ is multiplication.

\end{enumerate}

\item {\it comultiplication modules}: these are $R$-modules $M$ such that every $L \leq_R M$ is of the form $L= (0:_M I)$ for some ideal
$I$ of $R$, equivalently, $L = (0:_M (0:_R L))$. Examples of comultiplication modules include (by \cite{A-TF2007}):

\begin{enumerate}

\item the Pr\"{u}fer group $\mathbb{Z}(p^{\infty})$ for any prime $p$;

\item any semisimple commutative ring $R$ (considered canonically as an $R$-module);

\item any cocyclic module over a complete Noetherian local commutative ring.

\end{enumerate}

\end{enumerate}

\end{punto}

\begin{notation}
With $\mathcal{L}(M)$ ($\mathcal{L}^{\mathrm{f.i.}}(M)$) we denote the
lattice of (fully invariant) $R$-submodules of $M.$ Moreover, for every $%
L\leq _{R}M$ we set
\begin{equation*}
\begin{tabular}{lll}
$\mathcal{U}^{_{\mathrm{f.i.}}}(L)$ & $:=$ & $\{\widetilde{L}\leq _{R}M\mid
\widetilde{L}\supseteq L$ and $\widetilde{L}\leq _{R}^{\mathrm{f.i.}}M\};$
\\
$\mathcal{Q}^{_{\mathrm{f.i.}}}(L)$ & $:=$ & $\{\widetilde{L}\leq _{R}M\mid
\widetilde{L}\supseteq L$ and $\widetilde{L}/L\leq _{R}^{\mathrm{f.i.}}M/L\}.
$%
\end{tabular}
\end{equation*}
\qquad
\end{notation}

\begin{lemma}\label{ww-fi} \emph{(\cite[Lemma 17]{RRRF-AS2005})}
Let $L\leq _{R}^{\mathrm{f.i.}}M.$ Then $\mathcal{Q}^{_{\mathrm{f.i.}}}(L)\subseteq \mathcal{U}^{_{\mathrm{f.i.}}}(L)$
with equality in case $_{R}M$ is self-projective. In particular, if $M$ is self-projective and duo,
then $M/L$ is self-projective and duo.
\end{lemma}

\begin{punto}
By $\mathrm{Max}(M)$ ($\mathrm{Max}^{\mathrm{f.i.}}(M)$), we denote
the {\it possibly empty} class of maximal $R$-submodules of $M$ (the class of maximal $(R,S)$-subbimodules of $_R M_S$, which consists of all proper fully invariant $R$-submodules of $M$ which are not strictly contained in any
fully invariant proper $R$-submodule of $M$). For every $L\leq_{R}M,$ we set
\begin{equation*}
\begin{tabular}{lll}
$\mathcal{M}(L)$ & $:=$ & $\{K\in \mathrm{Max}(M)\mid K\supseteq L\};$ \\
$\mathcal{M}^{_{\mathrm{f.i.}}}(L)$ & $:=$ & $\{K\in \mathrm{Max}^{\mathrm{%
f.i.}}(M)\mid K\supseteq L\}.$%
\end{tabular}
\end{equation*}
\end{punto}

\begin{punto}
Let $L\leq _{R}M.$ We say that $L$ is \emph{essential} or \emph{large} in $M,$
and write $L\trianglelefteq M,$ iff $L\cap \widetilde{L}\neq 0$ for every $%
0\neq \widetilde{L}\leq M.$ On the other hand, we say $L$ is \emph{%
superfluous} or \emph{small} in $M,$ and we write $L\ll M$, iff $L+%
\widetilde{L}\neq M$ for every $\widetilde{L}\lvertneqq _{R}M.$ The \emph{%
radical} of $M$ is defined as
\begin{equation*}
\mathrm{Rad}(M):=\dbigcap\limits_{L\in \mathrm{Max}(M)}L=\dsum\limits_{L\ll
M}L,
\end{equation*}
whereas the \emph{socle} of $M$ is defined as
\begin{equation*}
\mathrm{Soc}(M):=\dsum\limits_{L\in \mathcal{S}(M)}L=\dbigcap\limits_{L%
\trianglelefteq M}L.
\end{equation*}
With $\mathrm{Spec}(R)$ we denote the spectrum of prime ideals
of $R.$
\end{punto}

\begin{definition}
We say $_{R} M$ is

\emph{local}, iff $M$ contains a proper $R$-submodule that contains every
proper $R$-submodule of $M$, equivalently iff $\dsum\limits_{L\lvertneqq
_{R}M}L\neq M$;

\emph{hollow} (or \emph{couniform}), iff for any $L_{1},L_{2}\lneqq _{R}M$ we have $L_{1}+L_{2}\lvertneqq _{R}M$, equivalently iff every proper $R$-submodule of $M$ is superfluous;

\emph{coatomic}, iff every proper $R$-submodule of $M$ is contained in a
maximal $R$-submodule of $M$, equivalently iff $\mathrm{Rad}(M/L)\neq M/L$
for every $L\lvertneqq _{R}M$;

\emph{f.i.-coatomic}, iff $\mathcal{M}^{\mathrm{f.i.}}(L)\neq \varnothing $
for every $L\lvertneqq _{R}^{\mathrm{f.i.}}M$, equivalently, iff $_R M_S$ is coatomic.
\end{definition}

\begin{exs} Finitely generated modules and semisimple modules are coatomic by \cite{Zos1980}. Hollow modules and semiperfect modules are coatomic by \cite{Gon1998}. If $R$ is left perfect (e.g. right Artinian), then every left $R$-module is coatomic by \cite{Gon1998}. The Pr\"{u}fer group $\mathbb{Z}_{p^{\infty }}$ is hollow, whence coatomic, but not local \cite[41.24 (6)]{Wis1991}.
\end{exs}

\begin{definition}
We call $_{R} M$ an {\it S-PCD}-module, iff $_R M$ is self-projective, coatomic and duo.
\end{definition}

\begin{lemma} Let $R$ be commutative, $M$ a finitely generated $R$-module and $\overline{R} := R/ann_R (M)$. The following are equivalent:

\begin{enumerate}

\item $_R M$ is an S-PCD-module;

\item $_R M$ is self-projective and duo;

\item $_{\overline {R}} M$ is projective and $_R M$ is duo;

\item $_{\overline {R}} M$ is projective and $_R M$ is locally cyclic;

\item $_R M$ is a multiplication module.

\end{enumerate}
\end{lemma}

\begin{Beweis}
The result follows directly from \cite[Theorem A]{Smi1994} noting that a finitely generated $R$-module $M$ is coatomic (e.g. \cite{Zos1980}) and that, moreover, $_R M$ is self-projective if and only if $_{\overline{R}} M$ is projective (e.g. \cite[Lemma 2.2]{AJK}).$\blacksquare$

\end{Beweis}

\begin{exs}\label{fgm}The following are examples of S-PCD-modules:

\begin{enumerate}
\item finitely generated multiplication modules over commutative rings (e.g. finitely generated self-projective ideals \cite[Proposition 9]{Smi1994}).

\item cyclic modules over left duo rings (such modules are quasi-projective by \cite[Theorem 1.2]{Koe1970} and are multiplication modules by \cite[Lemma 2.1]{Tug2003}).

\item finitely generated projective left ideals over a ring with commutative multiplication of left ideals (e.g. strongly regular rings);

\item left duo rings;

\item self-projective duo left modules over left perfect rings;

\item self-projective duo modules that are finitely generated (resp. semisimple, hollow, semiperfect).

\end{enumerate}
\end{exs}

\section{Fully Prime Submodules}

\qquad In this section, we recall, investigate and topologize the spectrum
of fully prime submodules of $M.$ For the special case $M=R$
a commutative ring, we recover several results on the classical Zariski
topology on the spectrum $\mathrm{Spec}(R)$ of prime ideals of $R$ (e.g.
\cite[pages 12-15]{AM1969}, \cite[II.4.3]{Bou1998}).

\begin{punto}
For any $X,Y\leq _{R}M$ consider
\begin{equation*}
X\ast _{M}Y:=\dsum_{f\in \mathrm{Hom}_{R}(M,Y)} f(X).
\end{equation*}
Notice that, if $Y\leq _{R}M$ is \emph{fully invariant}, then $X\ast
_{M}Y\leq _{R}M$ is also fully invariant; and if $X\leq _{R}M$ is fully
invariant, then $X\ast _{M}Y\subseteq X\cap Y.$
\end{punto}

\begin{definition}
We call $K\lvertneqq _{R}^{\mathrm{f.i.}}M$ a \emph{fully prime submodule}, or {\it fully prime in $M$}, iff for any $X,Y\leq _{R}^{\mathrm{f.i.}}M:$
\begin{equation*}
X\ast _{M}Y\subseteq K\Rightarrow X\subseteq K\text{ or }Y\subseteq K.
\end{equation*}
We call $_{R}M$ a \emph{fully prime} \emph{module}, iff $0$ is fully prime in $M$; equivalently, iff for any $X,Y\leq
_{R}^{\mathrm{f.i.}}M:$%
\begin{equation*}
X\ast _{M}Y=0\Rightarrow X=0 \text { or } Y=0.
\end{equation*}
\end{definition}

\begin{proposition}
\label{fp-cog}\emph{(\cite[1.6.3]{Wij2006}, \cite[3.1]{WW2009})} The
following are equivalent:

\begin{enumerate}
\item  $_{R}M$ is fully prime;

\item  $M$ is $K$-cogenerated for every $0 \neq K\leq _{R}^{\mathrm{f.i.}} M;$

\item  ${\rm Cogen}(M) = {\rm Cogen}(K)$ for any $0 \neq K\leq
_{R}^{\mathrm{f.i.}}M.$
\end{enumerate}
\end{proposition}

\begin{definition}
We call $_R M$ a {\it prime module}, iff $ann_R (M) = ann_R (L)$ for every $0 \neq L \leq_R M$. We call $K \lvertneqq_R M$ a {\it prime submodule}, or {\it prime in $M$}, iff $M/K$ is a prime $R$-module.
\end{definition}

The following result follows directly from the definition:

\begin{lemma}\label{K-prime} The following are equivalent for $K \lvertneqq_R M$:

\begin{enumerate}

\item $K$ is prime in $M$;

\item For any ideal $I$ of $R$ and $L \leq_R M$, we have: $$I L \subseteq K \Rightarrow L \subseteq K \text{  or  } I M \subseteq K;$$

\item For any $r \in R$ and $m \in M$, we have: $$rRm \subseteq K \Rightarrow m \in K \text{  or  } rM \subseteq K.$$
\end{enumerate}

\end{lemma}

\begin{lemma}\label{f.p->p}
If $K \lvertneqq_R^{\rm f.i.} M$ is fully prime in $M$, then $K$ is prime in $M$. In particular, if $_R M$ is fully prime, then $_R M$ is prime.
\end{lemma}

\begin{Beweis} Let $K \lvertneqq_R^{\rm f.i.} M$ be fully prime in $M$. Suppose that for some $r \in R$ and $m \in M$ we have $rRm \subseteq K$. Set $X := RrM$ and $Y := RmS$. Take any $x = \sum_{i=1}^n r_i r m_i \in X$. Let $f \in {\rm Hom}_R (M,Y)$ be arbitrary and suppose $f(m_i) = \sum_{j=1}^{k_i} r_{ij} m g_{ij} \in Y$. Then $$f(x) = f(\sum_{i=1}^n r_i r m_i) = \sum_{i=1}^n r_i r f(m_i) = \sum_{i=1}^n r_i r (\sum_{j=1}^{k_i} r_{ij} m g_{ij}) = \sum_{i=1}^n r_i g_{ij}(\sum_{j=1}^{k_i} r r_{ij} m) \in K.$$ It follows that $X *_M Y \subseteq K$. Since $K \lvertneqq_R^{\rm f.i.} M$ is fully prime in $M$, we conclude that $X \subseteq K$ whence $r M \subseteq K$, or $Y \subseteq K$ whence $m\in K$, i.e. $K$ is prime in $M$.$\blacksquare$
\end{Beweis}

\begin{remark}
The definition of fully prime modules we adopt is a modification of the definition of
prime modules in the sense of Bican et. al. \cite{BJKN80}, where arbitrary submodules are replaced by fully invariant ones. Following \cite{LP2000}, we call such modules B-\emph{prime}.  In fact, $_R M$ is B-prime if and only if $M$ is cogenerated by each of its non-zero $R$-submodules. Clearly, every B-prime module is fully prime. A duo module is B-prime if and only if it is fully prime. For more details on fully prime modules, the reader is referred to \cite{Wij2006} and \cite{WW2009}.
\end{remark}

\begin{exs}

\begin{enumerate}
\item Every homogeneous semisimple module is B-prime, whence fully prime.

\item $_{\mathbb{Z}}\mathbb{Q}$ is trivially fully prime, since it has no
non-trivial fully invariant $\mathbb{Z}$-submodules. However, $_{\mathbb{Z}}%
\mathbb{Q}$ is not B-prime, since it is not cogenerated by $\mathbb{Z}.$
\end{enumerate}
\end{exs}

\begin{punto}
Set
\begin{equation*}
\mathrm{Spec}^{\mathrm{fp}}(M):=\{K\lvertneqq _{R}^{\mathrm{f.i.}}M\mid K%
\text{ is a fully prime submodule of }M\}.
\end{equation*}

For every $L\leq _{R}M$ we set
\begin{equation*}
\mathcal{V}^{\mathrm{fp}}(L):=\{K\in \mathrm{Spec}^{\mathrm{fp}}(M)\mid
L\subseteq K\} \text {   ,  }  \mathcal{X}^{\mathrm{fp}}(L):=\{K\in \mathrm{Spec}^{\mathrm{fp}}(M)\mid
L \nsubseteqq K\}
\end{equation*}
and
\begin{equation*}
\mathrm{Rad}_M ^{\mathrm{fp}}(L):=\bigcap\limits_{K\in \mathcal{V}^{\mathrm{fp}%
}(L)}K\text{ \ \ (}:=M,\text{ iff }\mathcal{V}^{\mathrm{fp}}(L)=\varnothing
\text{).}
\end{equation*}
In particular, $\mathrm{Rad}_M ^{\mathrm{fp}}(M) = M$. Moreover, we set
\begin{equation*}
\mathrm{Rad} ^{\mathrm{fp}}(M):= \mathrm{Rad}_{M}^{\mathrm{fp}}(0) = \bigcap\limits_{K\in \mathrm{Spec}^{\mathrm{fp%
}}(M)}K\text{ }\;\;\text{(}:=M,\text{ iff }\mathrm{Spec}^{\mathrm{fp}%
}(M)=\varnothing \text{).}
\end{equation*}
We say that $L\leq _{R}^{\mathrm{f.i.}}M$ is $\mathrm{fp}$\emph{-radical}, iff
$\mathrm{Rad}_M ^{\mathrm{fp}}(L)=L.$ We say that $_{R}M$ is $\mathrm{fp}$\emph{-primeless}, iff $\mathrm{Spec}^{%
\mathrm{fp}}(M)=\varnothing$, equivalently iff $\mathrm{Rad} ^{\mathrm{fp}}(M)=M$.
\end{punto}

The following result can be derived from \cite[Proposition 18]{RRRF-AS2005}. We include its proof for convenience of the reader.

\begin{proposition}
\label{f.p.-pai}Let $L\lvertneqq _{R}^{\mathrm{f.i.}}M.$

\begin{enumerate}
\item  If $K\in \mathrm{Spec}^{\mathrm{fp}}(M)\cap \mathcal{Q}^{_{\mathrm{%
f.i.}}}(L),$ then $K/L\in \mathrm{Spec}^{\mathrm{fp}}(M/L).$

\item  If $_{R}M$ is self-projective, then the canonical epimorphism $M%
\overset{\pi }{\rightarrow }M/L\rightarrow 0$ induces bijections
\begin{equation*}
\mathrm{Spec}^{\mathrm{fp}}(M)\cap \mathcal{Q}^{_{\mathrm{f.i.}%
}}(L)\leftrightarrow \mathrm{Spec}^{\mathrm{fp}}(M/L).
\end{equation*}
\end{enumerate}
\end{proposition}

\begin{Beweis}
\begin{enumerate}
\item  Let $K\in \mathrm{Spec}^{\mathrm{fp}}(M)\cap \mathcal{Q}^{_{\mathrm{%
f.i.}}}(L),$ so that - by assumption - $K/L\leq _{R}^{\mathrm{f.i.}}M/L.$
Let $X/L,$ $Y/L\leq _{R}^{\mathrm{f.i.}}M/L$ be such that $X/L\ast
_{M/L}Y/L\subseteq K/L.$ Notice that $X,Y\in \mathcal{Q}^{_{\mathrm{f.i.}%
}}(L)\subseteq \mathcal{U}^{_{\mathrm{f.i.}}}(L)\ $by Lemma \ref{ww-fi}. For
every $f\in \mathrm{Hom}_{R}(M,Y),$ define
\begin{equation*}
\overline{f}:M/L\rightarrow Y/L,\text{ }\overline{m}\mapsto \overline{f(m)}.
\end{equation*}
Since $L\leq _{R}^{\mathrm{f.i.}}M,$ the map $\overline{f}$ is well-defined
and $\overline{f}(X/L)\subseteq K/L$ whence $f(X)\subseteq K$ and $X\ast _{M}Y\subseteq K.$ Since $K\in \mathrm{Spec}^{\mathrm{fp}}(M),$ $%
X\subseteq K$ so that $X/L\subseteq K/L$ or $Y\subseteq K$ so that $%
Y/L\subseteq K/L.$ Consequently, $K/L\in \mathrm{Spec}^{\mathrm{fp}}(M/L).$

\item  Assume that $_{R}M$ is self-projective. Let $K/L\in \mathrm{Spec}^{%
\mathrm{fp}}(M/L).$ Then, in particular, $K\in \mathcal{Q}^{_{\mathrm{f.i.}%
}}(L).$ Let $X,Y\leq _{R}^{\mathrm{f.i.}}M$ be such that $X\ast
_{M}Y\subseteq K.$ Since $_{R}M$ is self-projective, $X+L,$ $Y+L\in \mathcal{%
Q}^{_{\mathrm{f.i.}}}(L)$ by Lemma \ref{ww-fi}. Moreover, $_{R}M$ is $Y$%
-projective (e.g. \cite[18.2]{Wis1991}) and so there exists for every $g\in
\mathrm{Hom}_{R}(M/L,(Y+L)/L)$ some $\widetilde{g}\in \mathrm{%
Hom}_{R}(M,Y)$ such that $g(\overline{m})=\overline{\widetilde{g}(m)}$ for
every $m\in M.$ In particular, for every $x\in X$ and $l\in L$ we have $g(%
\overline{x+l})=\overline{\widetilde{g}(x+l)}\in K/L.$ Thus $(X+L)/L$ $\ast
_{M/L}$ $(Y+L)/L\subseteq K/L$ and it follows, since $K/L\in \mathrm{Spec}^{%
\mathrm{fp}}(M/L),$ that $(X+L)/L\subseteq K/L$ so that $X\subseteq K$ or $%
(Y+L)/L\subseteq K/L$ so that $Y\subseteq K.$ Consequently, $K\in \mathrm{%
Spec}^{\mathrm{fp}}(M).\blacksquare $
\end{enumerate}
\end{Beweis}

\qquad As a consequence of Proposition \ref{f.p.-pai}, we recover
\cite[1.6.3]{Wij2006}:

\begin{corollary}
\label{MKp}Let $K\lvertneqq _{R}^{\mathrm{f.i.}}M.$

\begin{enumerate}
\item  If $K$ is fully prime in $M,$ then $M/K$ is a fully prime $R$-module;

\item  If $_{R}M$ is self-projective, then
\begin{equation*}
K\text{ is fully prime in }M\Leftrightarrow M/K\text{ is a fully prime }R%
\text{-module.}
\end{equation*}
\end{enumerate}
\end{corollary}

\begin{remark}
\label{duo-coatomic}Notice that for every $L\in \mathrm{Max}^{\mathrm{f.i.}%
}(M),$ the $R$-module $M/L$ is trivially fully prime. If $_{R}M$ is
self-projective (and duo), then $\mathrm{Max}^{\mathrm{f.i.}}(M)\subseteq \mathrm{Spec}%
^{\mathrm{fp}}(M)$ ($\mathrm{Max}(M)\subseteq \mathrm{Spec}^{\mathrm{fp}}(M)$) by Proposition \ref{f.p.-pai}.
\end{remark}

Before we proceed, we introduce a class of modules which will prove to be useful in the sequel:

\begin{notation}
For any $L\in
\mathrm{Max}(M)$ set
\begin{equation}
L^{e}:=\dbigcap\limits_{\substack{ K\in \mathrm{Max}(M) \\ K\neq L}}K\text{
\ \ (}:=M\text{, iff }\mathrm{Max}(M)=\{L\}\text{).}  \label{Le}
\end{equation}
\end{notation}

\begin{punto}\label{max-prop}
We say that $_{R}M$ has the \emph{complete max-property}, iff for any $L\in
\mathrm{Max}(M)$ we have $L^{e}\nsubseteqq L.$ We also say that $_{R}M$ has the \emph{max-property}, iff for any  $L\in
\mathrm{Max}(M)$ and any finite subset $\mathcal{A} \in
\mathrm{Max}(M)\setminus\{L\}$ we have $\bigcap_{K\in \mathcal{A}} K \nsubseteqq L.$ Indeed, a module with a finite number of maximal submodules has the complete max-property if and only if it has the max property. Notice that $_{R}M$ satisfies the (complete) max-property if and only if $\mathrm{Max}(M)$ is {\it {\emph (}completely{\emph )} coindependent} in the sense of \cite[page 8]{CLVW2006}. A ring is said to have the (complete) max-property, iff its spectrum of maximal two-sided ideals is ({\it completely}) {\it coindependent}. For a survey on modules with the (complete) max-property, see \cite{Smi}.
\end{punto}

\begin{exs}
Every $R$-module with at most one maximal submodule (e.g. a local module)
has the complete max-property. On the other hand, each $R$-module with $L^{e} = 0$ for some $L \in {\rm Max}(M)$
does \emph{not} have the complete max-property. In particular, the ring of integers $\mathbb{Z}$
(considered as a $\mathbb{Z}$-module) does \emph{not} have the complete max-property.
\end{exs}

\begin{lemma}
\label{semi-local}Let $_{R}M$ be self-projective and duo with a finite number of maximal $R$-submodules. Then $M$ has the {\emph (}complete{\emph )} max-property.
\end{lemma}

\begin{Beweis}
If ${\rm Max}(M) = \varnothing$ we are done. Let $\mathrm{Max}(M)=\{M_{1},\cdots ,M_{n}\}.$ Since $_R M$ is duo and
self-projective, $\mathrm{Max}(M)\subseteq \mathrm{Spec}^{\mathrm{fp}}(M)$
(see Remark \ref{duo-coatomic}) and so $M_{i}^{e}$ $\nsubseteqq $ $M_{i}$
for each $i=1,\cdots ,n$ (otherwise $M_{j}=M_{i}$ for some $j\neq i$, a
contradiction).$\blacksquare $
\end{Beweis}

\begin{ex}(\cite[Corollary 3.8]{Smi}) Let $R$ be a commutative ring, $\{P_{\lambda}\}_{\lambda \in \Lambda}$ a non-empty collection of distinct maximal ideals of $R$ and $\{n_{\lambda}\}_{\lambda \in \Lambda}$ a collection of positive integers. Then $M:= \bigoplus_{\lambda \in \Lambda} R/{P^{n_{\lambda}}_{\lambda}}$ has the complete max-property.
\end{ex}

\begin{definition}
Let $K\in \mathrm{Spec}^{\mathrm{fp}}(M).$ We say that $K$ is \emph{minimal
above }$L,$ where $L\lvertneqq _{R}^{\mathrm{f.i.}}M,$ iff $K$ is a minimal
element of
\begin{equation*}
\mathcal{V}^{\mathrm{fp}}(L):=\{K\in \mathrm{Spec}^{\mathrm{fp}}(M)\mid
L\subseteq K\},
\end{equation*}
equivalently, iff $K$ contains $L$ and there is no $\widetilde{K}\in \mathrm{%
Spec}^{\mathrm{fp}}(M)$ that contains $L$ and is strictly contained in $K$.
We say that $K$ is \emph{minimal in }$\mathrm{Spec}^{\mathrm{fp}}(M),$ iff $K
$ is minimal above $0.$
\end{definition}

\begin{lemma}
\label{minimal}Let $_R M$ be self-projective and f.i.-coatomic. For every $%
L\lvertneqq _{R}^{\mathrm{f.i.}}M$ there exists $K\in \mathrm{Spec}^{\mathrm{%
fp}}(M)$ which is minimal above $L.$ In particular, $\mathrm{Spec}^{\mathrm{%
fp}}(M)$ has minimal elements.
\end{lemma}

\begin{Beweis}
Let $L\lvertneqq _{R}^{\mathrm{f.i.}}M.$ Since $_R M$ is self-projective and f.i.-coatomic, $%
\varnothing \neq \mathcal{M}^{\mathrm{f.i.}}(L)\subseteq \mathcal{V}^{%
\mathrm{fp}}(L).$ Let
\begin{equation*}
K_{1}\supseteq K_{2}\supseteq \cdots \supseteq K_{n}\supseteq
K_{n+1}\supseteq \cdots
\end{equation*}
be a descending chain in $\mathcal{V}^{\mathrm{fp}}(L)$ and set $%
K:=\bigcap\limits_{i=1}^{\infty }K_{i}.$ Suppose there exist $%
L_{1},L_{2}\leq _{R}^{\mathrm{f.i.}}M$ with $L_{1}\ast _{M}L_{2}\subseteq K$
but $L_{1}\nsubseteqq K$ and $L_{2}\nsubseteqq K.$ Then $L_{1}\nsubseteqq
K_{n_{1}}$ for some $n_{1}$ and $L_{2}\nsubseteqq K_{n_{2}}$ for some $n_{2}.
$ Setting $n:=\max \{n_{1},n_{2}\},$ we have $L_{1}\ast _{M}L_{2}\subseteq
K_{n}$ while $L_{1}\nsubseteqq K_{n}$ and $L_{2}\nsubseteqq K_{n}$ (a
contradiction). Therefore, $K\in \mathcal{V}^{\mathrm{fp}}(L).$ By Zorn's Lemma, $%
\mathcal{V}^{\mathrm{fp}}(L)$ has a minimal element. Applying this argument
to $L = 0$, we conclude that $\mathrm{Spec}^{\mathrm{fp}}(M) = \mathcal{V}^{\mathrm{fp}}(0)$ has a minimal element.$\blacksquare $
\end{Beweis}

\subsection*{Top$^{\mathrm{fp}}$-modules}

\begin{notation}
Set
\begin{equation*}
\begin{tabular}{cccccc}
$\xi ^{\mathrm{fp}}(M)$ & $:=$ & $\{\mathcal{V}^{\mathrm{fp}}(L)\mid L\in
\mathcal{L}(M)\};$ & $\xi _{\mathrm{f.i.}}^{\mathrm{fp}}(M)$ & $:=$ & $\{%
\mathcal{V}^{\mathrm{fp}}(L)\mid L\in \mathcal{L}^{\mathrm{f.i.}}(M)\};$ \\
$\tau ^{\mathrm{fp}}(M)$ & $:=$ & $\{\mathcal{X}^{\mathrm{fp}}(L)\mid L\in
\mathcal{L}(M)\};$ & $\tau _{\mathrm{f.i.}}^{\mathrm{fp}}(M)$ & $:=$ & $\{%
\mathcal{X}^{\mathrm{fp}}(L)\mid L\in \mathcal{L}^{\mathrm{f.i.}}(M)\};$ \\
$\mathbf{Z}^{\mathrm{fp}}(M)$ & $:=$ & $(\mathrm{Spec}^{\mathrm{fp}}(M),\tau
^{\mathrm{fp}}(M));$ & $\mathbf{Z}_{\mathrm{f.i.}}^{\mathrm{fp}}(M)$ & $:=$
& $(\mathrm{Spec}^{\mathrm{fp}}(M),\tau _{\mathrm{f.i.}}^{\mathrm{fp}}(M)).$%
\end{tabular}
\end{equation*}
\end{notation}

For an arbitrary $R$-module $M$, the set $\xi ^{\mathrm{fp}}(M)$ is not necessarily closed under finite
unions. This inspires the following:

\begin{definition}
We call $_R M$ a \emph{top}$^{\mathrm{fp}}$\emph{-module}, iff $\xi ^{\mathrm{fp%
}}(M)$ is closed under finite unions.
\end{definition}

\begin{ex}
It follows directly from the definition that uniserial modules are top$^{\mathrm{fp}}$-modules.
\end{ex}

\begin{lemma}
\label{fp-Properties}

\begin{enumerate}
\item  $\mathcal{V}^{\mathrm{fp}}(M)=\varnothing $ and $\mathcal{V}^{\mathrm{%
fp}}(0)=\mathrm{Spec}^{\mathrm{fp}}(M).$

\item  If $\{L_{\lambda }\}_{\Lambda }\subseteq \mathcal{L}(M),$ then $%
\bigcap\limits_{\Lambda }\mathcal{V}^{\mathrm{fp}}(L_{\lambda })=\mathcal{V}%
^{\mathrm{fp}}(\sum\limits_{\Lambda }L_{\lambda }).$

\item  If $L,\widetilde{L}\in \mathcal{L}^{\mathrm{f.i.}}(M),$ then
\begin{equation*}
\mathcal{V}^{\mathrm{fp}}(L)\cup \mathcal{V}^{\mathrm{fp}}(\widetilde{L})=%
\mathcal{V}^{\mathrm{fp}}(L\cap \widetilde{L})=\mathcal{V}^{\mathrm{fp}%
}(L\ast _{M}\widetilde{L}).
\end{equation*}
\end{enumerate}
\end{lemma}

\begin{Beweis}
Statements ``1'', ``2''\ and the inclusion $\mathcal{V}^{\mathrm{fp}}(L)\cup
\mathcal{V}^{\mathrm{fp}}(\widetilde{L})\subseteq \mathcal{V}^{\mathrm{fp}%
}(L\cap \widetilde{L})\subseteq \mathcal{V}^{\mathrm{fp}}(L\ast _{M}%
\widetilde{L})$ in ``3''\ are obvious. Conversely, if $K\in \mathcal{V}^{%
\mathrm{fp}}(L\ast _{M}\widetilde{L})$ then $L\ast _{M}\widetilde{L}%
\subseteq K.$ Since $K$ is fully prime in $M,$ we have $L\subseteq K$ so
that $K\in \mathcal{V}^{\mathrm{fp}}(L)$ or $\widetilde{L}\subseteq K$ so
that $K\in \mathcal{V}^{\mathrm{fp}}(\widetilde{L}).\blacksquare $
\end{Beweis}

\begin{theorem}
$\mathbf{Z}_{\mathrm{f.i.}}^{\mathrm{fp}}(M):=(\mathrm{Spec}^{\mathrm{fp}%
}(M),\tau _{\mathrm{f.i.}}^{\mathrm{fp}}(M))$ is a topological space. If $_R M$
is duo, then $M$ is a top$^{\mathrm{fp}}$-module and $\mathbf{Z}^{\mathrm{fp%
}}(M):=(\mathrm{Spec}^{\mathrm{fp}}(M),\tau ^{\mathrm{fp}}(M))$ is a
topological space.
\end{theorem}

\begin{notation}
For $\mathcal{A}\subseteq \mathrm{Spec}^{\mathrm{fp}}(M)$ set
\begin{equation*}
{\mathcal J}(\mathcal{A}):=\dbigcap_{K\in \mathcal{A}}K\text{ \ \ (}:=M,\text{
if }\mathcal{A}=\varnothing \text{).}
\end{equation*}
\end{notation}

\begin{lemma}\label{fp-closure} Let $M$ be a ${\rm top}^{\rm fp}$-module. The closure of any subset $\mathcal{A}\subseteq \mathrm{Spec}^{\mathrm{fp}%
}(M)$ is
\begin{equation}
\overline{\mathcal{A}}=\mathcal{V}^{\mathrm{fp}}({\mathcal J}(\mathcal{A})).  \label{fpclosure}
\end{equation}
\end{lemma}

\begin{Beweis} Let $\mathcal{A}\subseteq \mathrm{Spec}^{\mathrm{fp}}(M)$ and denote its closure by $\overline{%
\mathcal{A}}$. Since $\mathcal{A}%
\subseteq \mathcal{V}^{\mathrm{fp}}({\mathcal J}(\mathcal{A}))$ and $\mathcal{V}%
^{\mathrm{fp}}({\mathcal J}(\mathcal{A})$ is closed, we have $\overline{%
\mathcal{A}}\subseteq \mathcal{V}^{\mathrm{fp}}({\mathcal J}(\mathcal{A})).$
On the other hand, suppose $H\in \mathcal{V}^{\mathrm{fp}}({\mathcal J}(\mathcal{A}))\backslash \mathcal{A}$ and let $\mathcal{X}^{\mathrm{fp}}(L)$ be a neighborhood of $H,$ so that $H\nsupseteqq L.$ Then there exists $W\in
\mathcal{A}$ with $W\nsupseteqq L$ (otherwise $H\supseteq {\mathcal J}(\mathcal{A})\supseteq L,$ a contradiction), i.e. $W\in \mathcal{X}^{\mathrm{fp}%
}(L)\cap (\mathcal{A}\backslash \{H\})$ is a cluster point of $\mathcal{A}.$
Consequently, $\overline{\mathcal{A}}=\mathcal{V}^{\mathrm{fp}}({\mathcal J}(\mathcal{A})).\blacksquare $
\end{Beweis}

\begin{remarks}\label{fp-rms}
Let $M$ be a ${\rm top}^{\rm fp}$-module and consider the Zariski topology $\mathbf{Z}^{\mathrm{fp}}(M):=(\mathrm{Spec}^{%
\mathrm{fp}}(M),\tau ^{\mathrm{fp}}(M)).$

\begin{enumerate}
\item  $\mathbf{Z}^{\mathrm{fp}}(M)$ is a $T_{0}$ (Kolmogorov) space.

\item Setting $\mathcal{X}_{m}^{\mathrm{fp}}:=\mathcal{X}^{\mathrm{fp}}(Rm)$ for each $m\in M,$ the
set
\begin{equation*}
\mathcal{B}:=\{\mathcal{X}_{m}^{\mathrm{fp}}\mid m\in M\}
\end{equation*}
is a basis of open sets for the Zariski topology $\mathbf{Z}^{\mathrm{fp}}(M):$ Let $\mathcal{X}^{\mathrm{fp}}(L)$ be an open set in $\mathbf{Z}^{\mathrm{fp}}(M)$ and let $K\in \mathcal{X}^{\mathrm{fp}}(L)$. Then there exists some $m \in L \setminus\ K$, whence $K \in \mathcal{X}_{m}^{\mathrm{fp}} \subseteq \mathcal{X}^{\mathrm{fp}}(L)$.

\item If $L\in \mathrm{Spec}^{\mathrm{fp}}(M),$ then $\overline{\{L\}}= \mathcal{V}^{\mathrm{fp}}({\mathcal J}(\{L\})) =
\mathcal{V}^{\mathrm{fp}}(L).$ In particular, for any $K\in \mathrm{Spec}^{%
\mathrm{fp}}(M):$%
\begin{equation*}
K\in \overline{\{L\}}\Leftrightarrow K\supseteq L.
\end{equation*}

\item For any $L \leq_R M$, we have $L \subseteq \mathrm{Rad}_M ^{\mathrm{fp}}(L)$. Moreover, for any $L_1 \leq L_2 \leq _{R}M$, we have
\begin{equation}\label{rad-sub}
\mathrm{Rad}_M ^{\mathrm{fp}}(L_1) \subseteq \mathrm{Rad}_M ^{\mathrm{fp}}(L_2).
\end{equation}

Notice that, if ${\rm Spec}^{\rm fp} (M) \neq \varnothing$, then $\mathrm{Rad}_M ^{\mathrm{fp}}(M)= M \supsetneqq \mathrm{Rad} ^{\mathrm{fp}}(M)$.

\item For any $L\leq _{R}M$ we have
\begin{equation}\label{radrad}
\mathrm{Rad}_M ^{\mathrm{fp}}(\mathrm{Rad}_M ^{\mathrm{fp}}(L))=\mathrm{Rad}_M^{%
\mathrm{fp}}(L).
\end{equation}

\item If $_R M$ is self-projective, then $\mathcal{M}^{\mathrm{f.i.}}(L)\subseteq \mathcal{V}^{\mathrm{fp}}(L)$ for every $L\leq _{R}^{\mathrm{f.i.}}M$.

\item If $_R M$ is an S-PCD-module, then for every $L \leq_R M$ we have

\begin{enumerate}

\item $\mathcal{V}^{\mathrm{fp}}(L)=\varnothing$ if and
only if $L=M.$

\item  If $\mathcal{X}^{\mathrm{fp}}(L)=\varnothing ,$ then $L\subseteq
\mathrm{Rad}(M).$

\end{enumerate}

\item  Let $M\overset{\theta }{\simeq }N$ be an isomorphism of $R$-modules.
Then we have bijections
\begin{equation*}
\mathrm{Spec}^{\mathrm{fp}}(M)\longleftrightarrow \mathrm{Spec}^{\mathrm{fp}%
}(N).
\end{equation*}
In particular, we have $\theta (\mathrm{Rad} ^{\mathrm{fp}}(M))=\mathrm{Rad}^{%
\mathrm{fp}}(N).$ Moreover, we have a \emph{homeomorphism} $\mathbf{Z}^{%
\mathrm{fp}}(M)\approx \mathbf{Z}^{\mathrm{fp}}(N).$

\end{enumerate}
\end{remarks}



\begin{definition}
We call a topological space $X$ ({\it countably}) {\it compact}, iff every open cover of $X$ has a finite subcover. Countably compact spaces are also called {\it Lindelof spaces}. Note that some authors (e.g. \cite{Bou1966, Bou1998}) assume that compact spaces are in addition Hausdorff.
\end{definition}

\begin{punto}
A topological space ${\mathbf X}$ is said to be {\it Noetherian}, iff every ascending (descending) chain of open (closed) is stationary.
\end{punto}

\begin{notation}Let $M$ be a ${\rm top}^{\rm fp}$-module. We set
\begin{eqnarray*}
\mathcal{R}^{\mathrm{fp}}(M):=\{L\in \mathcal{L}(M)\mid \mathrm{Rad}_M^{%
\mathrm{fp}}(L)=L\};\\
\mathbf{CL}(\mathbf{Z}^{\mathrm{fp}}(M)):=\{\mathcal{A}\subseteq \mathrm{Spec%
}^{\mathrm{fp}}(M)\mid \mathcal{A}=\overline{\mathcal{A}}\}.
\end{eqnarray*}
\end{notation}

\begin{theorem}\label{noeth}Let $M$ be a ${\rm top}^{\rm fp}$-module.

\begin{enumerate}
\item We have a bijection between the class of \emph{fp}-radical $R$-submodules of $M$ and the class of closed sets in $\mathbf{Z}^{\mathrm{fp}}(M)$:
\begin{equation}\label{1-1}
\mathcal{V}^{\mathrm{fp}}(-) : \mathcal{R}^{\mathrm{fp}}(M) \longleftrightarrow \mathbf{CL}(%
\mathbf{Z}^{\mathrm{fp}}(M)),\text{ }L\mapsto \mathcal{V}^{\mathrm{fp}}(L).
\end{equation}

\item $\mathbf{Z}^{\mathrm{fp}}(M)$ is Noetherian if and only if $M$ has the ACC on its fp-radical submodules.

\item If $_R M$ is Noetherian, then
$\mathrm{Spec}^{\mathrm{fp}}(M)$ is Noetherian.
\end{enumerate}
\end{theorem}

\begin{Beweis}

\begin{enumerate}

\item Define
\begin{equation*}
\psi :\mathbf{CL}(\mathbf{Z}^{\mathrm{fp}}(M))\rightarrow \mathcal{R}^{%
\mathrm{fp}}(M),\text{ }\mathcal{V}^{\mathrm{fp}}(L)\mapsto \mathrm{Rad}_M^{%
\mathrm{fp}}(L).
\end{equation*}
For every $L\in \mathcal{R}^{\mathrm{fp}}(M)$ we have
\begin{equation*}
\psi (\mathcal{V}^{\mathrm{fp}}(L))=\mathrm{Rad}_M ^{\mathrm{fp}}(L)=L.
\end{equation*}
On the other hand, for every $\mathcal{A}=\mathcal{V}^{%
\mathrm{fp}}(K)\in \mathbf{CL}(\mathbf{Z}^{\mathrm{fp}}(M))$ we have
\begin{equation*}
\mathcal{V}^{\mathrm{fp}}(\psi (\mathcal{A}))= \mathcal{V}^{\mathrm{fp}} (\psi(\mathcal{V}^{%
\mathrm{fp}}(K)))=\mathcal{V}^{\mathrm{fp}}(\mathrm{Rad}_M ^{\mathrm{fp}}(K))=\mathcal{V}^{%
\mathrm{fp}}({\mathcal J}(\mathcal{A}))=\overline{\mathcal{A}}=\mathcal{A}%
.\blacksquare
\end{equation*}

\item This follows directly from "1".

\item This follows directly from "2". However, we provide here a direct proof. Assume that $_R M$ is Noetherian. Suppose that $\mathcal{X} (L)$ is an open set and consider an open basic cover $\{ \mathcal{X} (Rm_{\lambda}) | \, \lambda \in \Lambda  \}$ for $\mathcal{X} (L)$, so that $\mathcal{X} (L) \subseteq \bigcup_{\lambda \in \Lambda} \mathcal{X} (Rm_\lambda) = \mathcal{X} (\sum_{\lambda \in \Lambda} Rm_\lambda)$. Since $_R M$ is Noetherian, $N :=  \sum_{\lambda \in \Lambda} Rm_\lambda$ is finitely generated and so  there exists $\{m_{\lambda _1},  \cdots , m_{\lambda _n}\} \subseteq \{m_\lambda | \lambda \in \Lambda \}$ such that $N = \sum_i^n Rm_{\lambda _i}$. Clearly, $\{ \mathcal{X} (Rm_{\lambda _i}) | i =1 , \cdots n \}$ is a finite open cover for $\mathcal{X} (L)$, i.e. $\mathcal{X} (L)$ is compact. Consequently, $\mathrm{Spec}^{\mathrm{fp}}(M)$ is Noetherian by \cite[II.4.2, Proposition 9]{Bou1998}.$\blacksquare$
\end{enumerate}
\end{Beweis}

\begin{definition} (\cite{Bou1966}, \cite{Bou1998}) A non-empty topological space $\mathbf{X}$ is said to be
\begin{enumerate}
\item \emph{ultraconnected}, iff the intersection of any two non-empty closed subsets is non-empty.

\item \emph{irreducible} (or \emph{hyperconnected}), iff $\mathbf{X}$ is not the union of two proper closed subsets; equivalently, iff the
intersection of any two non-empty open subsets is non-empty.

\item \emph{connected}, iff $\mathbf{X}$ is not the {\it disjoint} union of two proper closed subsets; equivalently, iff the only subsets of $\mathbf{X}$ that are open and closed are $\varnothing $ and $\mathbf{X}$.
\end{enumerate}
\end{definition}

\begin{punto} (\cite{Bou1966}, \cite{Bou1998}) Let $\mathbf{X}$ be a non-empty topological space. A non-empty subset ${\mathcal A} \subseteq {\mathbf X}$ is an {\it irreducible set} in $\mathbf{X}$, iff it's an irreducible space w.r.t. the relative (subspace) topology; equivalently, iff for any proper closed subsets $A_1, A_2$ of $\mathbf X$ we have $${\mathcal A} \subseteq A_1 \cup A_2 \Rightarrow {\mathcal  A} \subseteq A_1 \text{  or  } {\mathcal A} \subseteq A_2.$$ A maximal irreducible subspace of $\mathbf{X}$ is called an \emph{irreducible component}. An irreducible component of a topological space is necessarily closed. The irreducible components of a Hausdorff space are just the singleton sets.
\end{punto}

\begin{proposition}\label{A-irred} Let $_R M$ be duo. Then $\mathcal{A} \subseteq \mathrm{Spec}^{\mathrm{fc}}(M)$ is irreducible if and only if ${\mathcal J}(\mathcal{A})$ is fully prime in $M$.
\end{proposition}

\begin{Beweis} Assume that $_R M$ is duo and let $\mathcal{A} \subseteq \mathrm{Spec}^{\mathrm{fc}}(M)$.

($\Rightarrow$): Assume that $\mathcal{A}$ is irreducible, whence - by definition - non-empty and consequently ${\mathcal J}(\mathcal{A}) \lvertneqq_R M$. Suppose that ${\mathcal J}(\mathcal{A})$ is not fully prime in $M.$ Then there exist $X,Y\leq_{R}M$ with ${\mathcal J}(\mathcal{A}) \supseteq X\ast _{M}Y$ but ${\mathcal J}(\mathcal{A})\nsupseteq X$ and ${\mathcal J}(\mathcal{A}) \nsupseteq Y.$ It follows that $\mathcal{A} \subseteq \mathcal{V}^{\mathrm{fp}}(X\ast _{M}Y)=\mathcal{V}^{\mathrm{fp}}(X)\cup \mathcal{V}^{\mathrm{fp}}(Y)$, a contradiction.
We conclude that ${\mathcal J}(\mathcal{A})$ is fully prime in $M.$

($\Leftarrow$): assume that ${\mathcal J}(\mathcal{A})$ is fully prime in $M$, whence - by definition - a proper $R$-submodule of $M$ and consequently $\mathcal{A} \neq \varnothing$. Suppose that $\mathcal{A} \subseteq \mathcal{V}^{\mathrm{fp}}(L_{1})\cup \mathcal{V}^{\mathrm{fp}}(L_{2})=\mathcal{V}^{\mathrm{fp}}(L_{1}\ast _{M}L_{2})$ for some $L_{1},L_{2}\leq _{R}M,$ so that
$L_{1}\ast _{M}L_{2}\subseteq {\mathcal J}(\mathcal{A}).$ Since, ${\mathcal J}(\mathcal{A})$ is fully prime in $M$, we conclude that ${\mathcal J}(\mathcal{A}) \supseteq L_{1},$ so that $\mathcal{A} \subseteq \mathcal{V}^{\mathrm{fp}}(L_{1})$; or ${\mathcal J}(\mathcal{A}) \supseteq L_{2},$ so that $\mathcal{A} \subseteq \mathcal{V}^{\mathrm{fp}}(L_{2})$. We conclude that $\mathcal{A}$ is irreducible.$\blacksquare $
\end{Beweis}

\begin{corollary}\label{rad-irred} Let $_R M$ be duo.

\begin{enumerate}

\item $\mathrm{Spec}^{\mathrm{fp}}(M)$ is irreducible if and only if $\mathrm{Rad} ^{\mathrm{fp}}(M)$ is fully prime in $M$.

\item If $_R M$ is self-projective, then ${\rm Max}(M)$ is irreducible if and only if ${\rm Rad}(M)$ is fully prime in $M$.
\end{enumerate}
\end{corollary}

\begin{exs} Let $_R M$ be duo.

\begin{enumerate}
\item Let $X \subseteq \mathrm{Spec}^{\mathrm{fp}}(M)$ is a chain, then $X$ is irreducible.

\item If $_R M$ is uniserial, then $\mathrm{Spec}^{\mathrm{fp}}(M)$ is irreducible.

\end{enumerate}
\end{exs}



\begin{proposition}\label{K-irred} Let $_{R} M$ be duo. The bijection \emph{(\ref{1-1})} restricts to bijections

\begin{equation}
\mathrm{Spec}^{\mathrm{fp}}(M) \longleftrightarrow \{ {\mathcal Y} | \, {\mathcal Y} \subseteq \mathrm{Spec}^{\mathrm{fp}}(M) \text{  is an irreducible closed subset} \}
\end{equation}

and

\begin{equation}
{\rm Min}(\mathrm{Spec}^{\mathrm{fp}}(M)) \longleftrightarrow \{ {\mathcal Y} | \, {\mathcal Y} \subseteq \mathrm{Spec}^{\mathrm{fp}}(M) \text{  is an irreducible component} \}.
\end{equation}

\end{proposition}

\begin{Beweis}

Recall from (\ref{1-1}) the bijection $\mathcal{R}^{\mathrm{fp}}(M) \stackrel{\mathcal{V}^{\mathrm{f.p.}}(-)}{\longleftrightarrow } \mathbf{CL}(\mathbf{Z}^{\mathrm{fp}}(M))$. If $K\in \mathrm{Spec}^{\mathrm{fp}}(M)$, then ${\mathcal J}({\mathcal{V}^{\rm fp}(K)}) = K$, and so $\mathcal{V}^{\rm fp}(K)$ is an irreducible closed set by Proposition \ref{A-irred}. Conversely, let ${\mathcal Y} = \mathcal{V}^{\mathrm{fp}}(L)$ be an irreducible closed set in $\mathbf{Z}^{\mathrm{fp}}(M).$ Then ${\mathcal Y} = \overline{\mathcal Y} = \mathcal{V}^{\mathrm{fp}}({\mathcal J}({\mathcal Y}))$ and so ${\mathcal J} ({\mathcal Y}) = {\rm Rad}^{\rm fp}_M (L)$ is fully prime in $M$ by Proposition \ref{A-irred}.

On the other hand, notice that $\mathrm{Spec}^{\mathrm{fp}}(M)$ has minimal elements by
Lemma \ref{minimal}. If $K$ is minimal in $\mathrm{Spec}^{\mathrm{fp}}(M)$,
then clearly $\mathcal{V}^{\mathrm{fp}}(K)$ is an irreducible component of $%
\mathrm{Spec}^{\mathrm{fp}}(M)$ by the argument above. Conversely, let $\mathcal{Y}$ be
an irreducible component of $\mathrm{Spec}^{\mathrm{fp}}(M).$ Then $\mathcal{%
Y}$ is closed and so $\mathcal{Y}=\mathcal{V}^{\mathrm{fp}}(L)$ for some $%
L\leq _{R}M.$ Since $\mathcal{Y}$ is irreducible, $K:= {\rm Rad}^{\rm fp}(L)\in \mathrm{Spec}^{%
\mathrm{fp}}(M)$ by ``1''. If $K$ is not minimal in $\mathrm{Spec}^{\mathrm{%
fp}}(M),$ then there exists $\widetilde{K}\in \mathrm{Spec}^{\mathrm{fp}}(M)$
such that $\widetilde{K}\subsetneqq K.$ In this case $\mathcal{V}^{\mathrm{fp%
}}(L) = \mathcal{V}^{\mathrm{fp}}(K) \subsetneqq \mathcal{V}^{\mathrm{fp}}(\widetilde{K})$ (a contradiction,
since $\mathcal{V}^{\mathrm{fp}}(\widetilde{K})$ is irreducible).$\blacksquare $
\end{Beweis}

\begin{definition} Let ${\mathbf X}$ be a topological space and $\mathcal{Y} \subset {\mathbf X}$ be an closed set. A point $y \in {\mathcal Y}$ is said to be a {\it generic point}, iff $ {\mathcal Y} = \overline{ \{y\}}$. If every irreducible closed subset of ${\mathbf X}$ has a unique generic point, then we call ${\mathbf X}$ a {\it Sober} space.
\end{definition}

\begin{corollary}Let $_R M$ be duo. Then $\mathrm{Spec}^{\mathrm{fp}}(M)$ is a Sober space.
\end{corollary}

\begin{proposition}\label{ultra} Let $M$ be an S-PCD-module. Then $_{R} M$ is hollow if and only if $\mathrm{Spec}^{\mathrm{fp}}(M)$ is ultraconnected.
\end{proposition}

\begin{Beweis}
Assume that $_{R} M$ is hollow. If $\mathcal{V}^{\mathrm{fp}}(L_{1}),$ $%
\mathcal{V}^{\mathrm{fp}}(L_{2})\subseteq \mathrm{Spec}^{\mathrm{fp}}(M)$
are any two non-empty closed subsets, then $L_{1}\neq M\neq L_{2}.$ Since $M$
is hollow, $L_{1}+L_{2}\neq M$ and so there exists $N\in \mathrm{Max}(M)$
such that $L_{1}+L_{2}\subseteq $ $N\subsetneqq M.$ Indeed, $\mathcal{V}^{%
\mathrm{fp}}(L_{1})\cap \mathcal{V}^{\mathrm{fp}}(L_{2})=\mathcal{V}^{%
\mathrm{fp}}(L_{1}+L_{2})\neq \varnothing ,$ since $N\in \mathcal{V}^{%
\mathrm{fp}}(L_{1}+L_{2}).$ Conversely, assume that the intersection of any
two non-empty closed subsets of $\mathrm{Spec}^{\mathrm{fp}}(M)$ is
non-empty. Let $L_{1},L_{2}\lvertneqq _{R}M,$ so that $\mathcal{V}^{\mathrm{%
fp}}(L_{1})\neq \varnothing \neq \mathcal{V}^{\mathrm{fp}}(L_{2}).$ By
assumption $\mathcal{V}^{\mathrm{fp}}(L_{1}+L_{2})=\mathcal{V}^{\mathrm{fp}%
}(L_{1})\cap \mathcal{V}^{\mathrm{fp}}(L_{2})\neq \varnothing ,$ whence $L_{1}+L_{2}\neq M$.
Consequently, $_{R} M$ is hollow.$\blacksquare $
\end{Beweis}

\begin{lemma}\label{open-compact} Let $_R M$ be duo.

\begin{enumerate}

\item Every countably compact open set of $\mathrm{Spec}^{\mathrm{fp}}(M)$ is of the form ${\mathcal X}(N)$ for some countably generated $R$-submodule $N\leq_R M$.

\item Every compact open set of $\mathrm{Spec}^{\mathrm{fp}}(M)$ is of the form ${\mathcal X}(N)$ for some finitely generated $R$-submodule $N\leq_R M$.

\end{enumerate}

\end{lemma}

\begin{Beweis}

\begin{enumerate}

\item Let $L \leq_R M$ an suppose that $\mathcal{X}^{\rm fp} (L)$ is a countably compact open set. Notice that $\mathcal{X}^{\rm fp} (L) = \mathcal{X}^{\rm fp} (\sum_{l \in L} Rl) = \bigcup_{l \in L} \mathcal{X}^{\rm fp} (Rl)$. Since $\mathcal{X}^{\rm fp} (L)$ is countably compact,  there exists a countable subset $\{l_i\}_{i=1}^\infty \subseteq L$ such that
$\mathcal{X}^{\rm fp} (L) = \bigcup_{l \in L} \mathcal{X}^{\rm fp} (Rl) = \mathcal{X}^{\rm fp} (\sum_{i=1}^\infty Rl_i)$.

\item The proof is analogous to that of "1".$\blacksquare$
\end{enumerate}
\end{Beweis}

\begin{theorem}
\label{fp-Lindelof}
Let $M$ be an S-PCD-module.

\begin{enumerate}

\item If $\mathrm{Max}(M)$ is countable, then $\mathbf{Z}^{\mathrm{fp}}(M)$ is countably compact.

\item If $\mathrm{Max}(M)$ is finite, then $\mathbf{Z}^{\mathrm{fp}}(M)$ is compact.
\end{enumerate}
\end{theorem}

\begin{Beweis} We need only to prove (1), since (2) can be proved analogously.

Assume that $\mathrm{Max}(M)=\{N_{\lambda _{k}}\}_{k\geq 1}$ is countable.
Let $\{\mathcal{X}^{\mathrm{fp}}(L_{\alpha })\}_{\alpha \in I}$ be
an open cover of $\mathrm{Spec}^{\mathrm{fp}}(M)$ (i.e. $\mathrm{Spec}^{%
\mathrm{fp}}(M)\subseteq \dbigcup\limits_{\alpha \in I}\mathcal{X}^{\mathrm{%
fp}}(L_{\alpha })$). Since $\mathrm{Max}(M)\subseteq \mathrm{Spec}^{\mathrm{%
fp}}(M)$ we can pick for each $k\geq 1,$ some $\alpha _{k}\in I$ such that $%
N_{\lambda _{k}}\nsupseteqq L_{\alpha _{k}}.$ Suppose $\dsum\limits_{k\geq
1}L_{\alpha _{k}}\neq M.$ Then there exists some $N\in \mathrm{Max}(M)$ such
that $M\neq N\supseteq \dsum\limits_{k\geq 1}L_{\alpha _{k}}$ (a
contradiction, since $N=N_{\lambda _{k}}\nsupseteq L_{\alpha _{k}}$ for some
$k\geq 1$). Hence $\dsum\limits_{k\geq 1}L_{\alpha _{k}}=M$ and we conclude
that $\mathrm{Spec}^{\mathrm{fp}}(M)=\mathcal{X}^{\mathrm{fp}}(M)=\mathcal{X}%
^{\mathrm{fp}}(\dsum\limits_{k\geq 1}L_{\alpha _{k}})=\dbigcup\limits_{k\geq
1}\mathcal{X}^{\mathrm{fp}}(L_{\alpha _{k}})$, i.e. $\{\mathcal{X}^{\mathrm{%
fp}}(L_{\alpha _{k}})\mid k\geq 1\}\subseteq \{\mathcal{X}^{\mathrm{fp}%
}(L_{\alpha })\}_{\alpha \in I}$ is a countable subcover.$\blacksquare $
\end{Beweis}

\begin{proposition}\label{local---conn} Let $_R M$ be duo and assume that $\mathrm{Spec}^{\mathrm{fp}}(M)=\mathrm{Max}(M).$

\begin{enumerate}

\item If $_R M$ has the complete max-property, then $\mathrm{Spec}^{\mathrm{fp}}(M)$ is discrete.

\item $M$ has a unique maximal submodule if and only if $M$ has the complete max-property and $\mathrm{Spec}^{\mathrm{fp}}(M)$ is connected.

\end{enumerate}
\end{proposition}

\begin{Beweis}
\begin{enumerate}

\item If $_R M$ has the complete max-property, then for every $K \in \mathrm{Spec}^{\mathrm{fp}}(M)$ we have $\{K\}=\mathcal{X}^{\mathrm{fp}}(\{K\}^{e})$, whence open. Consequently, $\mathrm{Spec}^{\mathrm{fp}}(M)$ is discrete.

\item If $_R M$ has a unique maximal submodule, it has indeed the complete max-property and since $\mathrm{Spec}^{\mathrm{fp}}(M)$ has only one point, it's indeed connected. On the other hand, if $M$ has the complete max-property, then $\mathbf{Z}^{\mathrm{fp}}(M)$ is discrete by (1) and so $|{\rm Max}(M)| = |\mathrm{Spec}^{\mathrm{fp}}(M)| = 1$ since a connected discrete space has only one point.$\blacksquare $
\end{enumerate}
\end{Beweis}

\begin{corollary}\label{local}
Let $_R M$ be an S-PCD-module and assume that every fully prime $R$-submodule of $M$ is maximal.

\begin{enumerate}

\item If $_R M$ has the complete max-property, then $\mathrm{Max}(M)$ is countable if and only if $\mathrm{Spec}^{\mathrm{fp}}(M)$ is countably compact.

\item $\mathrm{Max}(M)$ is finite if and only if ${_R M}$ has the complete max-property and $\mathrm{Spec}^{\mathrm{fp}}(M)$ is compact.

\item $_R M$ is local if and only if $_R M$ has the max property and $\mathrm{Spec}^{\mathrm{fp}}(M)$ is connected.

\end{enumerate}
\end{corollary}

\begin{lemma}Let $_R M$ be an S-PCD-module. If $n\geq 2$ and $\mathcal{A}=\{K_{1},...,K_{n}\}\subseteq \mathrm{Spec}^{%
\mathrm{fp}}(M)$ is connected, then for every $i\in \{1,...,n\},$
there exists $j\in \{1,...,n\}\backslash \{i\}$ such that $K_{i}\leq
_{R}K_{j}$ or $K_{j}\leq _{R}K_{i}.$
\end{lemma}

\begin{Beweis}
Without loss of generality, suppose $K_{1}\nsubseteqq K_{j}$ and $%
K_{j}\nsubseteqq K_{1}$ for all $2\leq j\leq n$ and set $F:=\dbigcap%
\limits_{i=2}^{n}K_{i},$ $W_{1}:=\mathcal{A}\cap \mathcal{X}^{\mathrm{fp}%
}(K_{1})=\{K_{2},...,K_{n}\}$ and $W_{2}:=\mathcal{A}\cap \mathcal{X}^{%
\mathrm{fp}}(F)=\{K_{1}\}$ (if $n=2,$ then clearly $W_{2}=\{K_{1}\};$ if $n>2
$ and $K_{1}\notin W_{2},$ then $(K_{2}\ast
_{M}\dbigcap\limits_{i=3}^{n}K_{i})\subseteq
\dbigcap\limits_{i=2}^{n}K_{i}\subseteq K_{1}$ and it follows that $\dbigcap\limits_{i=3}^{n}K_{i}\subseteq K_{1}.$ One can show
by induction that $K_{n}\subseteq K_{1},$ a contradiction). So $\mathcal{A}%
=W_{1}\cup W_{2},$ a disjoint union of proper non-empty open subsets, a
contradiction.$\blacksquare $
\end{Beweis}

\begin{definition}
A collection $\mathcal{G}$ of subsets of a topological space $\mathbf{X}$ is
\emph{locally finite}, iff every point of $\mathbf{X}$ has a neighborhood
that intersects only finitely many elements of $\mathcal{G}.$
\end{definition}

\begin{proposition}
Let $M$ be an S-PCD-module and have the complete max-property. Let $\mathcal{K}=\{K_{\lambda
}\}_{\Lambda }\subseteq \mathrm{Max}(M)$ be non-empty. If $\left| \mathcal{M}%
(L)\right| <\infty $ for every $L\in \mathrm{Spec}^{\mathrm{fp}}(M),$ then $%
\mathcal{K}$ is locally finite.
\end{proposition}

\begin{Beweis}
Let $L\in \mathrm{Spec}^{\mathrm{fp}}(M)$ and set
\begin{equation*}
F:=\bigcap\limits_{K\in \mathcal{K}\cap \mathcal{X}^{\mathrm{fp}}(L)}K\text{
\ \ (}:=M,\text{ iff }\mathcal{K}\cap \mathcal{X}^{\mathrm{fp}%
}(L)=\varnothing \text{)}.
\end{equation*}
Notice that $F\nsubseteqq L:$ If $F\subseteq L,$ then there exists a
maximal $R$-submodule $F\subseteq L\subseteq \widetilde{K}\subsetneqq M.$
Since $M$ has the complete max-property, we conclude that $\widetilde{%
K}=K$ for some $K\in \mathcal{K}\cap \mathcal{X}^{\mathrm{fp}}(L)$, a
contradiction. Therefore, $L\in \mathcal{X}^{\mathrm{fp}}(F).$ It follows directly
from the assumptions and the construction of $F$ that
$\mathcal{K}\cap \mathcal{X}^{\mathrm{fp}}(F) = \mathcal{K}\cap \mathcal{V}^{\mathrm{fp}}(L) \subseteq  \mathcal{M}(L)$, whence finite.$\blacksquare$
\end{Beweis}

\begin{lemma}If $_R M$\ is an S-PCD-module, then the following are equivalent for any $L \leq_{R} M$:
\begin{enumerate}

\item $L\in \mathrm{Max}(M)$;

\item $L$ is fully prime in $M$ and $\mathcal{V}^{\mathrm{fp}}(L)=\{L\}$;

\item $\{L\}$ is a closed set in $\mathbf{Z}_{M}^{\mathrm{fp}}.$

\end{enumerate}
\end{lemma}

\begin{Beweis}
Recall first that $\mathrm{Spec}^{\mathrm{fp}}(M) \supseteq {\rm Max}(M)$ by Remark \ref{duo-coatomic}. Let $L \leq_{R} M$.

(a) $\Rightarrow$ (b) and (b) $\Rightarrow$ (c) are obvious.

(c) $\Rightarrow$ (a): Assume that $\{L\}$ is closed in $\mathbf{Z}^{\mathrm{fp}}(M)$, so that $\{L\} = \mathcal{V}^{\mathrm{fp}}(K)$ for some $K \leq_{R} M$. If $L \notin {\rm Max}(M)$, then $L \subsetneqq \tilde{L}$ for some $\tilde{L} \in {\rm Max}(M)$. In this case, $\{L, \tilde{L}\} \subseteq \mathcal{V}^{\mathrm{fp}}(K) = \{L\}$, a contradiction. Consequently, $L \in \mathrm{Max}(M)$.$\blacksquare$
\end{Beweis}

The following result follows directly from the pervious lemma and the fact that a topological space is  $T_{1}$ if and only if every singleton subset is closed.

\begin{proposition}\label{Frecht} If $_R M$\ is an S-PCD-module, then $\mathrm{Spec}^{\mathrm{fp}}(M)=\mathrm{Max}(M)$ if and only if $\mathbf{Z}^{\mathrm{fp}}(M)$ is $T_{1}$ \emph{(}Fr\'{e}cht space\emph{)}.
\end{proposition}

Combining the assertions in Propositions \ref{local---conn} and \ref{Frecht}, we obtain:

\begin{theorem}\label{f.p.-discrete}
If $_R M$\ is an S-PCD-module and has the complete max-property, then the following are equivalent:
\begin{enumerate}
\item  $\mathrm{Spec}^{\mathrm{fp}}(M)=\mathrm{Max}(M);$

\item  $\mathbf{Z}^{\mathrm{fp}}(M)$ is discrete;

\item  $\mathbf{Z}^{\mathrm{fp}}(M)$ is $T_{2}$ \emph{(}Hausdorff space\emph{%
)};

\item  If $\mathbf{Z}^{\mathrm{fp}}(M)$ is $T_{1}$ \emph{(}Fr\'{e}cht space%
\emph{)}.
\end{enumerate}
\end{theorem}

\subsection*{${\textbf{Top}}^{\rm fp}$ Rings}

\qquad In what follows, we give some applications to the Zariski topology on the prime spectrum of an associative (not necessarily commutative) ring. We write $_R R$ ($R_R$) to indicate that we consider $R$ as a left (right) $R$-module. Notice that $_R R$ ($R_R$) is self-projective and coatomic, whence an S-PCD-module if and only if $R$ is left (right) duo. To avoid repetition, we include only some results which distinguish $_R R$ from arbitrary modules.

\begin{punto} Consider the ring $R$ with the canonical left $R$-module structure.
Notice that we have an isomorphism of rings
\begin{equation*}
R\simeq \mathrm{End}(_{R}R)^{op} \text{   and   } R\simeq \mathrm{End}(R_R).
\end{equation*}
The submodules of $_{R}R$ ($R_R$) which are fully invariant in $R$ coincide with the
{\it two-sided} ideals of $R.$ For two ideals $I,J$ of $R,$ we have $I\ast
_{R}J=IJ$ (the usual product of ideals). Therefore, ${\rm Spec}^{\rm fp}(_R R) = {\rm Spec}(R) = {\rm Spec}^{\rm fp}(R_R)$, where ${\rm Spec}(R)$ is the spectrum of prime {\it two-sided} ideals of $R$.
In particular, $_{R}R$ is fully prime if and only if $R$ is a prime ring if and only if $R_R$ is fully prime.
\end{punto}

Notice that ${\rm Spec}(R)$ attains a Zariski topology declaring the closed sets to be $$\xi ^{\mathrm{fp}}(R) = \{ {\mathcal V}(I) | \, I \text{ is an ideal of } R\}, \text{ where } {\mathcal V}(I) = \{J | J \text{ is an ideal of } R \text{ and } I \subseteq J \}.$$ However, $$\xi ^{\mathrm{fp}}(_R R) = \{ {\mathcal V}(I) | \, I \text{ is a left ideal of } R \} \text{ and } \xi ^{\mathrm{fp}}(R_R) = \{ {\mathcal V}(I) | \, I \text{ is a right ideal of } R \}$$ are not necessarily closed under finite unions.

\begin{definition}
We call the ring $R$ {\it left} ({\it right}) ${top}^{\rm fp}$-ring, iff $\xi ^{\mathrm{fp}}(_R R)$ ($\xi ^{\mathrm{fp}}(R_R)$) is closed under finite unions, equivalently iff $\mathbf{Z}^{\mathrm{fp}}(_R R)$ ($\mathbf{Z}^{\mathrm{fp}}(R_R)$) is a topological space. We call $R$ a {\it top}$^{\rm fp}$-{\it ring}, iff $R$ is both a left and a right ${\rm top}^{\rm fp}$-ring.
\end{definition}

\begin{exs}

\begin{enumerate}

\item If $R$ is left (right) uniserial, then $R$ is a left (right) ${\rm top}^{\rm fp}$-ring. If $R$ is uniserial, then $R$ is a ${\rm top}^{\rm fp}$-ring.

\item If $R$ is left (right) duo, then $R$ is a left (right) ${\rm top}^{\rm fp}$-ring. If $R$ is duo, then $R$ is a ${\rm top}^{\rm fp}$-ring.


\end{enumerate}
\end{exs}

\begin{proposition}
Let $R$ be a left {\rm (}right{\rm )} ${\rm top}^{\rm fp}$-ring. Then $\mathbf{Z}^{\mathrm{fp}}(_R R)$ $(\mathbf{Z}^{\mathrm{fp}}(R_R))$ is compact.
\end{proposition}

\begin{Beweis}
Assume that $R$ is a left ${\rm top}^{\rm fp}$-ring and let $\{ {\mathcal X}_{a_{\lambda}} | \lambda \in \Lambda \}$ be a basic open cover for ${\rm Spec}(R)$, so that ${\rm Spec}(R) = \bigcup_{\lambda \in \Lambda}  {\mathcal X}_{a_{\lambda}}$. Then $\varnothing = \bigcap_{\lambda \in \Lambda} {\mathcal V}(Ra_\lambda) =  {\mathcal V}(\sum_{\lambda \in \Lambda} Ra_\lambda)$, whence $\sum_{\lambda \in \Lambda} Ra_\lambda = R$. It follows that there exist $\{r_{\lambda_1}, \cdots, r_{\lambda_n}\} \subseteq R$ and $\{a_{\lambda_1}, \cdots, a_{\lambda_n} \}, \lambda_i \in \Lambda$ such that $\sum_{i=1}^n r_{\lambda_i} a_{\lambda_i} = 1$. Clearly, $\varnothing = \bigcap_{i=1}^n {\mathcal V}(Ra_i)$, and so $\{ {\mathcal X}_{a_{{\lambda}_i}} | i=1, \cdots n\}$ is a finite subcover.$\blacksquare$
\end{Beweis}

\begin{punto}
The ring $R$ is called $\pi$-{\it regular}, iff for each $a \in R$ there exist a positive
integer $n = n(a)$, depending on $a$, and $x \in R$ such that $a^n = a^n x a^n$. If $R$ is a left (right) duo ring, then every prime ideal of $R$ is maximal if and only if $R$ is $\pi$-regular \cite{Hir1978}.
\end{punto}

\begin{corollary}
Assume that $R$ is left {\rm (}right{\rm )} duo and $\pi$-regular. Then $R$ has a finite number of maximal ideals if and only if $R$ has the complete max-property.
\end{corollary}

\begin{corollary}\label{pi-reg}
Assume that $R$ is left {\rm (}right{\rm )} duo. The following are equivalent:

\begin{enumerate}

\item $R$ is $0$-dimensional and has a finite number of maximal ideals;

\item $R$ is $\pi$-regular and has the complete max property;

\item $\mathrm{Spec}(R)$ is discrete;

\item $\mathrm{Spec}(R)$ is $T_2$ and finite;

\item $\mathrm{Spec}(R)$ is $T_1$ and finite.

\end{enumerate}
\end{corollary}

We finish this section with an application  to commutative rings. Recall that a ring is {\it semilocal}, iff $R/{\rm Rad}(R)$ is semisimple. A commutative ring is semilocal if and only if it has a finite number of maximal ideals.

\begin{corollary}\label{0-dim}
Let $R$ be a commutative ring.

\begin{enumerate}

\item The following are equivalent:

\begin{enumerate}

\item $R$ is $0$-dimensional and semilocal;

\item $R$ is $\pi$-regular and has the complete max property;

\item $\mathrm{Spec}(R)$ is discrete;

\item $\mathrm{Spec}(R)$ is $T_2$ and finite;

\item $\mathrm{Spec}(R)$ is $T_1$ and finite.

\end{enumerate}

\item $R$ is semisimple if and only if $R$ is von Neumann regular and satisfies any of the following conditions:

\begin{enumerate}

\item $R$ is semilocal;

\item $R$ has the complete max-property;

\item $R$ is Noetherian;

\item $R$ is perfect.

\end{enumerate}
\end{enumerate}
\end{corollary}

\textbf{Acknowledgement}: The author thanks Professor Patrick Smith for fruitful discussions on the topic during his visit to the University of Glasgow (September 2008) and thereafter.

\end{document}